\documentclass[notitlepage,a4paper,fleqn,9pt]{extarticle}

\usepackage[T1]{fontenc}
\usepackage{fancyhdr}
\usepackage{indentfirst}
\usepackage{times}
\usepackage[dvips]{graphicx} 

\usepackage{amsmath}
\usepackage{amsfonts}

\title{\bf Numerical Solutions of a Boundary Value Problem for the Anomalous Diffusion Equation\\
with the Riesz Fractional Derivative }

\author{\bf Mariusz Ciesielski and Jacek Leszczynski}

\date{\vspace{-2ex}\em 
Institute of Mathematics and Computer Science, 
Czestochowa University of Technology \\
ul. Dabrowskiego 73, 42-200 Czestochowa \\
e-mail:  mariusz@imi.pcz.pl, jale@imi.pcz.pl } 

\fancypagestyle{plain}

\makeatletter
\renewcommand{\@seccntformat }[1]{\csname the#1\endcsname.\quad}
\renewcommand\section{\@startsection {section}{1}
{-\parindent}{3ex \@plus -1ex \@minus -.2ex}
{3ex \@plus -1ex \@minus -.2ex}{\textbf}}
\renewcommand\subsection{\@startsection {subsection}{2}
{-\parindent}{1.5ex \@plus -1ex \@minus -.2ex}
{1.5ex \@plus -1ex \@minus -.2ex}{\textit}}
\renewcommand\subsubsection{\@startsection {subsubsection}{2}
{-\parindent}{1.5ex \@plus -1ex \@minus -.2ex}
{1.5ex \@plus -1ex \@minus -.2ex}{\textit}} 
\makeatother

\begin{document}
\twocolumn[ 
\maketitle
\noindent\rule{\textwidth}{.5pt}

\noindent
Abstract\\
\vspace{1ex}

\noindent 
In this paper we present in one-dimensional space a numerical solution of a 
partial differential equation of fractional order. This equation describes a~process 
of anomalous diffusion. The process arises from the interactions within the complex 
and non-homogeneous background. We presented a numerical method which bases on the 
finite differences method. We considered pure initial and boundary-initial value problems 
for the equation with the Riesz-Feller fractional derivative. In the final part
of this paper sample results of simulation were shown. \\
\vspace{1ex}


\noindent\rule{\textwidth}{.5pt}
\vspace{2ex}
]{

\section{Introduction} 
Anomalous diffusion is a phenomenon strongly connected
with the interactions within complex and non-homoge\-neous background. This
phenomenon is observed in transport of fluid in porous materials, in the
chaotic heat baths, amorphous semiconductors, particle dynamics inside
polymer network, two-dimen\-sional rotating flow and also in econophysics.
Phenomenon of anomalous diffusion deviates from the standard diffusion behaviour.
In opposite to standard diffusion where linear form in the mean square 
displacement $\left\langle \,x^{2}\left( t\right) \right\rangle \sim k_{1}t $
of diffusing particle over time occurs,
anomalous diffusion is characterized by the non-linear one
$\left\langle \,x^{2}\left( t\right) \right\rangle \sim k_{\gamma }t^{\gamma} $,
for $\gamma \in (0,2]$.
In this phenomenon may exist dependence 
$\left\langle \,x^{2}\left( t\right) \right\rangle \rightarrow \infty $ ,
which is characterized by occurrence of rare but
extremely large jumps of diffusing particle -- well-known as the Levy motion
or the Levy flights. Ordinary diffusion follows Gaussian statistics and
Fick's second law for finding running process at time $t$ whereas anomalous
diffusion follows non-Gaussian statistic or can be interpreted as the Levy
stable densities.

Many authors proposed models which base on linear and non-linear forms
of differential equations. Such models can simulate anomalous diffusion but
they don't reflect its real behaviour. Several authors~\cite{Carpinteri,  
Gorenflo2, Gorenflo1, Hilfer1, Mainardi, Metzler1, Podlubny} apply
fractional calculus in modelling of this type of diffusion. This means that
time and spatial derivatives in the classical diffusion equation are
replaced by fractional ones. In comparison to derivatives of integer
order, which depend on the local behaviour of the function, derivatives of
fractional order accumulate the whole history of this function.

\section{Mathematical background}
In this paper, we consider an equation in the following form 
\setlength\arraycolsep{5pt}
\begin{equation}
\frac{\partial }{\partial t}C(x,t)=k_{\alpha }\frac{\partial ^{\alpha }}
{\partial \left\vert x\right\vert ^{\alpha }}C(x,t) \text{, } t\geq 0 \text{, } x\in \mathbb{R} \text{, }
\end{equation}
where 
$C(x,t)$ is a field variable, $\frac{\partial ^{\alpha }}
{\partial \left\vert x\right\vert ^{\alpha }}C(x,t)$ is the Riesz-Feller fractional
operator~\cite{Metzler1, Samko}, 
$\alpha $ is the real order of this operator, $k_{\alpha }$ is the
coefficient of generalized (anomalous) diffusion with the unit of measure 
$\left[ m^{\alpha }/s\right]$.
According to~\cite{Gorenflo1, Mainardi} the Riesz-Feller fractional operator 
for $0<\alpha \leq 2$, $\alpha \neq 1$ for one-variable function $u(x)$ is
\begin{equation}
\begin{split}
\frac{\partial ^{\alpha }}{\partial \left\vert x\right\vert ^{\alpha }}%
u(x)={}_{x}D_{\theta }^{\alpha }u\left( x\right) =& {}-\left[ c_{L}\left(
\alpha ,\theta \right) \,_{-\infty }D_{x}^{\alpha }u\left( x\right) \right. 
\\
& \quad \ \left. +c_{R}\left( \alpha ,\theta \right) \,_{x}D_{+\infty}^{\alpha }
u\left( x\right) \right] \hbox{,}
\end{split}
\end{equation}
where 
\begin{equation}
{}_{-\infty }D_{x}^{\alpha }u\left( x\right) =
\begin{cases}
\dfrac{d}{dx}\left[ _{-\infty }I_{x}^{1-\alpha }u\left( x\right) \right] , & %
\mbox{for}\ 0<\alpha \leq 1\text{,} \vspace{0.2cm}\\ 
\dfrac{d^{2}}{dx^{2}}\left[ _{-\infty }I_{x}^{2-\alpha }u\left( x\right) %
\right] , & \mbox{for}\ 1<\alpha \leq 2\mbox{,}%
\end{cases}%
\end{equation}%
\begin{equation}
{}_{x}D_{+\infty }^{\alpha }u\left( x\right) =%
\begin{cases}
-\dfrac{d}{dx}\left[ _{x}I_{+\infty }^{1-\alpha }u\left( x\right) \right] ,
& \mbox{for}\ 0<\alpha \leq 1\text{,} \vspace{0.2cm}\\ 
\dfrac{d^{2}}{dx^{2}}\left[ _{x}I_{+\infty }^{2-\alpha }u\left( x\right) %
\right] , & \mbox{for}\ 1<\alpha \leq 2\text{.}%
\end{cases}%
\end{equation}%
and coefficients $c_{L}\left( \alpha ,\theta \right) $, $c_{R}\left( \alpha
,\theta \right) $ (for $0<\alpha \leq 2$, $\alpha \neq 1$, and for $%
\left\vert \theta \right\vert \leq \min \left( \alpha ,2-\alpha \right) $),
are defined as
\begin{equation}
c_{L}\left( \alpha ,\theta \right) =\frac{\sin \dfrac{\left( \alpha -\theta
\right) \pi }{2}}{\sin \left( \alpha \pi \right) }\text{, }\quad c_{R}\left(
\alpha ,\theta \right) =\frac{\sin \dfrac{\left( \alpha +\theta \right) \pi 
}{2}}{\sin \left( \alpha \pi \right) }\text{.}  \label{coeff_c}
\end{equation}

The fractional operators of order $\alpha $: $_{-\infty}I_{x}^{\alpha }u\left( x\right) $
and $_{x}I_{\infty }^{\alpha }u\left(x\right) $
are defined as the left- and right-side of Weyl fractional integrals
~\cite{Gorenflo2, Gorenflo1, Oldham, Podlubny, Samko} which definitions are
\begin{equation}
_{-\infty }I_{x}^{\alpha }u\left( x\right) = \frac{1}{\Gamma (\alpha )}
\int_{-\infty }^{x} \frac{u\left( \xi \right) }{(x-\xi )^{1-\alpha }}d\xi 
\text{,}
\end{equation}
\begin{equation}
_{x}I_{\infty }^{\alpha }u\left( x\right) = \frac{1}{\Gamma(\alpha )}
\int_{x}^{\infty } \frac{u\left( \xi \right) }{(\xi-x)^{1-\alpha }}d\xi 
\text{.}
\end{equation}

\vspace{0.3cm}
Considering Eqn~(1) we obtain the classical diffusion equation for 
$\alpha =2$, i.e. the heat transfer equation. If $\alpha =1$, and the parameter of
skewness $\theta $ admits extreme values in~(\ref{coeff_c}), the transport
equation is noted. Therefore we assume variations of the parameter $\alpha $
within the range $0<\alpha \leq 2$. Analysing behaviour of the parameter 
$\alpha <2$ in Eqn~(1), we found some combination between transport and
propagation processes. 

For analytic solution of Eqn~(1) we can apply Green functions~\cite{Gorenflo2}.
We numerically solve Eqn~(1) when additional
non-linear term may occur. Some numerical methods used in solution of
fractional differential equations can be found in~\cite{Gorenflo1}. However
they apply the infinite domain. 

In this work we will consider Eqn (1) limited for $1 <\alpha \leq 2$ 
in one dimensional domain $\Omega :L\leq x\leq R$ with the
boundary-value conditions of the first kind (the Dirichlet conditions) as
\begin{equation}
\left\{ 
\begin{array}{ll}
x=L: & C\left( L,t\right) =g_{L}\left( t\right) , \vspace{0.1cm}\\ 
x=R: & C\left( R,t\right) =g_{R}\left( t\right) ,
\end{array}
\right. t>0,
\end{equation}
and with the initial-value condition
\begin{equation}
\left. C\left( x,t\right) \right\vert _{t=0}=c_{0}\left( x\right) \text{.}
\end{equation}

\section{Numerical method}
According to the finite difference method \cite{Ames, Frank1, Frank2, Hoffman, Majchrzak}
we consider a~discrete from of Eqn~(1)
both in time and space. In the previous work \cite{Ciesielski} we solved numerically 
the anomalous diffusion equation similar 
to the Eqn~(1) with the time-fractional derivative. We called this method FFDM 
(Fractional FDM). The problem of solving of Eqn (1) lies in properly approximation
of the Riesz-Feller derivative~(2) in numerical scheme. 

\subsection{Approximation of the Riesz-Feller derivative}
We begin numerical analysis from discrete forms of operators~(6) and~(7). 
We introduce homogenous spatial grid $-\infty <\ldots
<x_{i-2}<x_{i-1}<x_{i}<x_{i+1}<x_{i+2}<\ldots <\infty $ with the step 
$h=x_{k}-x_{k-1}$ and we denote value of function $u$ in the point $x_{k}$ as 
$u_{k}=u\left( x_{k}\right)$, for $k\in \mathbb{Z}$. 
In order to simplify notations we take here the function of one variable.
For numerical integration scheme we assumed the trapezoidal rule.
The integral~(6) in point $x_{i}$ of the grid is replaced by the sum of
discrete integrals as 
\begin{equation}
_{-\infty }I_{x_{i}}^{\alpha }u_{i}=\dfrac{1}{\Gamma (\alpha )}%
\sum\limits_{k=0}^{\infty }\int\limits_{x_{i-k-1}}^{x_{i-k}}\dfrac{u\left(
\xi \right) }{(x_{i}-\xi )^{1-\alpha }}d\xi \text{,}
\end{equation}
and using linear interpolation of function $u$ in every sub-interval $%
[x_{i-k-1},x_{i-k}]$ 
\begin{equation}
u^{\ast }\left( \xi \right) =\dfrac{u_{i-k}-u_{i-k-1}}{h}\xi +\dfrac{%
u_{i-k-1}x_{i-k}-u_{i-k}x_{i-k-1}}{h}
\end{equation}%
we have
\setlength\arraycolsep{1pt}
\begin{eqnarray}
_{-\infty }I_{x_{i}}^{\alpha }u_{i} &\approx &\dfrac{1}{\Gamma (\alpha )}
\sum\limits_{k=0}^{\infty }\int\limits_{x_{i-k-1}}^{x_{i-k}}\dfrac{u^{\ast
}\left( \xi \right) }{(x_{i}-\xi )^{1-\alpha }}d\xi  \notag \\
&=&
\dfrac{1}{\Gamma (\alpha )}  \sum\limits_{k=0}^{\infty } 
\left[ \left(u_{i-k}-u_{i-k-1}\right) a_{k}^{\left( \alpha \right) }\right. \\ 
&& \hspace{1.5cm} \left. +\left( u_{i-k-1}x_{i-k}-u_{i-k}x_{i-k-1}\right) b_{k}^{\left(\alpha \right) }\right] \notag
\end{eqnarray}
where 
\setlength\arraycolsep{1pt}
\begin{eqnarray}
a_{k}^{\left( \alpha \right) } &=&h^{\alpha -1}x_{i}\dfrac{\left( k+1\right)
^{\alpha }-k^{\alpha }}{\alpha }-h^{\alpha }\dfrac{\left( k+1\right)
^{1+\alpha }-k^{1+\alpha }}{1+\alpha }\text{,} \\
b_{k}^{\left( \alpha \right) } &=&h^{\alpha -1}\dfrac{\left( k+1\right)
^{\alpha }-k^{\alpha }}{\alpha }\text{.}
\end{eqnarray}

\vspace{2cm}
\noindent After next transforms we can write
\begin{equation}
_{-\infty }I_{x_{i}}^{\alpha }u_{i}\approx h^{\alpha
}\sum\limits_{k=0}^{\infty }u_{i-k}v_{k}^{\left( \alpha \right) }
\end{equation}%
where
\setlength\arraycolsep{0pt}
\begin{eqnarray}
&& v_{k}^{\left( \alpha \right) }= \dfrac{1}{\Gamma (2+\alpha )}  \times \\
&& \left\{ 
\begin{array}{ll}
1 \text{,} & \text{ for }k=0 \text{,} \vspace{0.2cm} \\ 
 \left( k+1\right) ^{1+\alpha }-2k^{1+\alpha }+\left( k-1\right)^{1+\alpha } \text{,} & \text{ for }k=1,...,\infty \text{.}
\end{array} \notag
\right. 
\end{eqnarray}

Similar to previous considerations we approximate operator 
$_{x}I_{\infty }^{\alpha }u\left( x\right) $
in the point $x_{i}$ and finally we obtain
\begin{equation}
_{x_{i}}I_{\infty }^{\alpha }u_{i}\approx h^{\alpha
}\sum\limits_{k=0}^{\infty }u_{i+k}v_{k}^{\left( \alpha \right) }\,\text{,}
\end{equation}
where coefficients $v_{k}^{\left( \alpha \right) }$ have identical forms as~(16). 

In the next step we analyse operator (2). It can be expressed in the form
(in order to simplify this we denote $c_{L}=c_{L}\left( \alpha ,\theta \right) $ 
and $c_{R}=c_{R}\left( \alpha ,\theta \right) $ )
\setlength\arraycolsep{0pt}
\begin{eqnarray}
&&_{x}D_{\theta }^{\alpha }u\left( x\right) ={} \\
&&-\left[ c_{L}\,\dfrac{d^{2}}{dx^{2}}\,\left[ _{-\infty }I_{x}^{2-\alpha
}u\left( x\right) \right] +{}{}c_{R}\,\dfrac{d^{2}}{dx^{2}}\,\left[
_{x}I_{+\infty }^{2-\alpha }u\left( x\right) \right] \right] \text{.}  \notag
\end{eqnarray}
We used the central difference scheme for the second spatial derivative 
in the point $x_{i}$ and we obtain
\setlength\arraycolsep{0pt}
\begin{eqnarray}
&&_{x_{i}}D_{\theta }^{\alpha }u_{i}\approx \\
&&
\begin{array}[t]{c}
-\left[ c_{L}\,\dfrac{_{-\infty }I_{x}^{2-\alpha }u_{i-1}-2\,_{-\infty
}I_{x}^{2-\alpha }u_{i}+\,_{-\infty }I_{x}^{2-\alpha }u_{i+1}}{h^{2}}\right.
\, \\ 
\left. +{}{}c_{R}\,\dfrac{_{x}I_{+\infty }^{2-\alpha
}u_{i-1}-2\,_{x}I_{+\infty }^{2-\alpha }u_{i}+\,_{x}I_{+\infty }^{2-\alpha
}u_{i+1}}{h^{2}}\right] \text{.}%
\end{array}
\notag
\end{eqnarray} 
After numerous transforms we obtain the final form as 
\begin{equation}
_{x_{i}}D_{\theta }^{\alpha }u_{i} \approx \dfrac{1}{%
h^{\alpha }}\sum\limits_{k=-\infty }^{\infty }u_{i+k}w_{k}^{\left( \alpha
\right) }\text{,}
\end{equation}%
where coefficients $w_{k}^{\left( \alpha \right) }$ are 
\setlength\arraycolsep{0pt}
\begin{eqnarray}
&&w_{k}^{\left( \alpha \right) }=\dfrac{-1}{\Gamma \left( 4-\alpha \right) }
\times \\
&&\left\{ 
\begin{array}{ll}
\left( \left( \left\vert k\right\vert +2\right) ^{3-\alpha }-4\left(
\left\vert k\right\vert +1\right) ^{3-\alpha }+6\left\vert k\right\vert
^{3-\alpha }\right. &  \\ 
\left. \quad \quad -4\left( \left\vert k\right\vert -1\right) ^{3-\alpha}
  +\left( \left\vert k\right\vert -2\right) ^{3-\alpha }\right) c_{L}, & \text{for }k\leq -2 \vspace{0.1cm}\\ 
\left( 3^{3-\alpha }-2^{5-\alpha }+6\right) c_{L}+c_{R}, & \text{for }k=-1 \vspace{0.1cm}\\ 
\left( 2^{3-\alpha }-4\right) \left( c_{L}+c_{R}\right), & \text{for }k=0 \vspace{0.1cm}\\ 
\left( 3^{3-\alpha }-2^{5-\alpha }+6\right) c_{R}+c_{L}, & \text{for }k=1 \vspace{0.1cm}\\ 
\left( \left( k+2\right) ^{3-\alpha }-4\left( k+1\right) ^{3-\alpha}+6k^{3-\alpha }\right. &  \\ 
\left. \quad \quad -4\left( k-1\right) ^{3-\alpha }+\left( k-2\right)
^{3-\alpha }\right) c_{R}, & \text{for }k\geq 2
\end{array}
\right. \text{.} \notag
\end{eqnarray}

Assuming $\alpha =2$ and $\theta =0$ we have
$c_{L}\left( 2,0\right) \,=$\ $c_{R}\left( 2,0\right) \,=-\frac{1}{2}$
and we obtain 
\setlength\arraycolsep{5pt}
\begin{equation}
w_{k}^{\left( 2\right) }=\left\{ 
\begin{array}{ll}
0, & \text{for }k\leq -2 \\ 
1, & \text{for }k=-1 \\ 
-2, & \text{for }k=0 \\ 
1, & \text{for }k=1 \\ 
0, & \text{for }k\geq 2%
\end{array}
\right. \text{.}
\end{equation}
These coeeficients are identical as for wide known the central difference scheme 
for the second derivative. Also when $\alpha \rightarrow 1^{+}$ and $\theta =0$
after arduous calculations of limits we obtain coefficients
\setlength\arraycolsep{0pt}
\begin{eqnarray}
&&w_{k}^{\left( 1^{+}\right) }=\dfrac{1}{2\pi }\times \\
&&\left\{ 
\begin{array}{ll}
\ln \dfrac{\left( \left\vert k\right\vert +1\right) ^{4\left( \left\vert
k\right\vert +1\right) ^{2}}\left( \left\vert k\right\vert -1\right)
^{4\left( \left\vert k\right\vert -1\right) ^{2}}}{\left( \left\vert
k\right\vert +2\right) ^{\left( \left\vert k\right\vert +2\right)
^{2}}\left\vert k\right\vert ^{6k^{2}}\left( \left\vert k\right\vert
-2\right) ^{\left( \left\vert k\right\vert -2\right) ^{2}}}, & \text{ for } k\leq -2, \\ 
16\ln 2-9\ln 3, & \text{ for }k=-1, \\ 
-8\ln 2, & \text{ for }k=0, \\ 
16\ln 2-9\ln 3, & \text{ for }k=1, \\ 
\ln \dfrac{\left( k+1\right) ^{4\left( k+1\right) ^{2}}\left( k-1\right)
^{4\left( k-1\right) ^{2}}}{\left( k+2\right) ^{\left( k+2\right)
^{2}}k^{6k^{2}}\left( k-2\right) ^{\left( k-2\right) ^{2}}}, & \text{ for }k\geq 2.%
\end{array}%
\right.  \notag
\end{eqnarray}
In literature didn't find exact values of approximating coefficients.
When $\alpha = 1$ the Riesz-Feller operator is singular, hence the problem.
Numerous works of Gorenflo and Mainardi i.e.~\cite{Gorenflo2, Gorenflo1}
propose various ways which determine values of the coefficients $w_{k}^{\left( \alpha \right) }$
(i.e. based on the Gr\"unwald-Letnikov discretization) but they don't
provide continuity in the interval $\alpha \in (1,2]$.
The coefficients~(23) can approximate the Cauchy process 
when we use~(23) in numerical calculations.

\subsection{Fractional FDM}
While discretization of the Riesz-Feller derivative in space is done,
in this subsection we describe the finite difference method for the
equation of anomalous diffusion (1). Here we restrict this solution to
one dimensional space. In comparison with the standard diffusion equation
where discretization of the second derivative in space can be
approximated by the central difference of second order, we will use
generalized scheme given by formula~(20). The differences appear in
setting of boundary conditions. 

We shall introduce a temporal grid $0=t^{0}<t^{1}<\ldots <t^{f}<t^{f+1}<\ldots <$$\infty $
with the step $\Delta t=t^{f+1}-t^{f}$ and we
denote value of the function $C\left( x,t\right) $ in the point $x_{k}$ at the moment of
time $t^{f}$ as $C_{k}^{f}=C\left( x_{k},t^{f}\right) $ for $k\in \mathrm{\!}
\mathbb{Z}$ and $f\in \mathrm{\mathbb{N}}$. 

\subsubsection{Pure initial value problem}
In the explicit scheme of the FDM we replaced Eqn~(1) by the following formula 
\begin{equation}
\dfrac{C_{i}^{f+1}-C_{i}^{f}}{\Delta t}=K_{\alpha }\dfrac{1}{h^{\alpha }}%
\sum\limits_{k=-\infty }^{\infty }C_{i+k}^{f}w_{k}^{\left( \alpha \right) }%
\text{.}
\end{equation}%
After simplification finally we obtained
\begin{equation}
C_{i}^{f+1}=\sum\limits_{k=-\infty }^{\infty }C_{i+k}^{f}p_{k}^{\left(
\alpha \right) }\text{,}
\end{equation}%
where coefficients $p_{k}^{\left( \alpha \right) }$ are 
\begin{equation}
p_{k}^{\left( \alpha \right) }=\left\{ 
\begin{array}{ll}
1+K_{\alpha }\dfrac{\Delta t}{h^{\alpha }}w_{0}^{\left( \alpha \right) }, & 
\text{for }k=0\text{,} \vspace{0.1cm} \\ 
K_{\alpha }\dfrac{\Delta t}{h^{\alpha }}w_{k}^{\left( \alpha \right) }, & 
\text{for }k\neq 0\,.%
\end{array}
\right.
\end{equation}

Using simple calculations one may proof, that arise the following relationship 
\begin{equation}
\sum_{k=-\infty }^{\infty }p_{k}^{\left( \alpha \right) }=1\text{.}
\end{equation}

In order to determine stability of the explicit scheme the
coefficient~(26) for $k=0$ in formula~(25) should be positive 
\begin{equation}
p_{0}^{\left( \alpha \right) }=1+K_{\alpha }\dfrac{\Delta t}{h^{\alpha }}%
w_{0}^{\left( \alpha \right) }>0\text{.}
\end{equation}
Hence we fixed the maximum length of the step $\Delta t$ as
\begin{equation}
\Delta t<\dfrac{-h^{\alpha }}{K_{\alpha }w_{0}^{\left( \alpha \right) }}=%
\dfrac{h^{\alpha }\Gamma \left( 4-\alpha \right) }{K_{\alpha }\left(
2^{3-\alpha }-4\right) \left( c_{L}\left( \alpha ,\theta \right)
+c_{R}\left( \alpha ,\theta \right) \right) }\text{.}
\end{equation}

The initial condition~(9) is introduced directly to every grid
nodes at the first step $t=t^{0}$. This determines values of the function~$C$ as
\begin{equation}
C_{i}^{0}=c_{0}\left( x_{i}\right) \text{.}
\end{equation}

In unbounded domains the implicit method isn't easily applicable because 
it generates infinite dimensions of all matrices. Thus one
usually seeks improved difference equations within the explicit scheme.

\subsubsection{Boundary-initial value problem}
Presenting numerical solution~(25) with included unbounded domain 
$-\infty <x<\infty $  has no practical implementations in computer simulations.

Now, we present solution of this problem on the finite domain 
$\Omega :L\leq x\leq R$ with boundary conditions~(8). We divide this
domain $\Omega $ into $N$ sub-domains with $h = (R-L)/N$. Figure~1 shows modified
spatial grid.
\begin{figure}[h]
\begin{center}
 \includegraphics[width=0.48\textwidth]{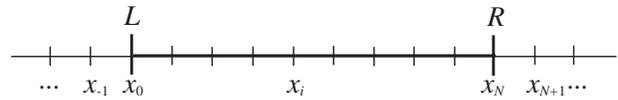}  
\end{center}
\vspace{-0.5cm}  
\caption{ The nodes grid over space}
\end{figure}

\noindent Here we can observe additional 'virtual' points in the grid
placed outside of the domain $\Omega $. In order to introduce the
Dirichlet boundary conditions we proposed treatment which bases
on assumption that values of the function $C$ in outside points are identical 
as values in the boundary nodes $x_{0}$ or $x_{N}$ 
\setlength\arraycolsep{5pt}
\begin{equation}
C\left( x_{k},t\right) =\left\{ 
\begin{array}{ll}
C\left( x_{0},t\right) =g_{L}\left( t\right) \text{,} & \text{for }k<0\text{,} \vspace{0.1cm} \\ 
C\left( x_{N},t\right) =g_{R}\left( t\right) \text{,} & \text{for }k>N\text{.}
\end{array}
\right. 
\end{equation}
On the base of previous considerations we modify expression~(20) for discretization of the
Riesz-Feller derivative. Thus we have
{\setlength\arraycolsep{1pt}
\begin{eqnarray}
_{x_{i}}D_{\theta }^{\alpha }C\left( x_{i},t\right)  &\approx &\dfrac{1}
{h^{\alpha }}\left[ \sum\limits_{k=-i}^{N-i}C\left( x_{i+k},t\right)
w_{k}^{\left( \alpha \right) }\right.   \notag \\
&&+\left. g_{L}\left( t\right) {s_{L}}_{i}^{\left( \alpha \right)
}+g_{R}\left( t\right) {s_{R}}_{i}^{\left( \alpha \right) }\right] ,
\end{eqnarray}}
for $i=1,\ldots ,N-1$, where
\setlength\arraycolsep{0pt}
\begin{eqnarray}
&&{s_{L}}_{i}^{\left( \alpha \right) }=\sum\limits_{k=-\infty}^{-i-1}
w_{k}^{\left( \alpha \right) }=\dfrac{-1}{\Gamma \left( 4-\alpha\right) }
\times  \\
&&\left[ -\left( i+2\right) ^{3-\alpha }+3\left( i+1\right) ^{3-\alpha}
-3i^{3-\alpha }+\left( i-1\right) ^{3-\alpha }\right] c_{L},  \notag
\end{eqnarray}
\begin{eqnarray}
&&{s_{R}}_{i}^{\left( \alpha \right) }=\sum\limits_{k=N-i+1}^{\infty}
w_{k}^{\left( \alpha \right) }=\dfrac{-1}{\Gamma \left( 4-\alpha \right) }
\left[ -\left( N-i+2\right) ^{3-\alpha }\right.  \\
&&\left. +3\left( N-i+1\right) ^{3-\alpha }-3\left( N-i\right) ^{3-\alpha}
+\left( N-i-1\right) ^{3-\alpha }\right] c_{R}.  \notag
\end{eqnarray}
Putting this expression to Eqn (1) we obtain a finite difference
scheme depending on weighting factor $\sigma $. Here we assumed
\setlength\arraycolsep{1pt}
\begin{eqnarray}
g_{L}^{f+\frac{1}{2}} &=& g_{L}\left( t^{f+\frac{1}{2}}\right) =g_{L}\left(
           \Delta t\left( f+\frac{1}{2}\right) \right), \\
g_{R}^{f+\frac{1}{2}} &=& g_{R}\left( t^{f+\frac{1}{2}}\right) 
=g_{R}\left( \Delta t\left( f+\frac{1}{2}\right) \right) 
\end{eqnarray}
in order to simplify the numerical scheme. For internal nodes $x_{i}$, $i=1,\ldots ,N-1$
we have
\setlength\arraycolsep{1pt}
\begin{eqnarray}
\dfrac{C_{i}^{f+1}-C_{i}^{f}}{\Delta t} &=&K_{\alpha }\dfrac{1}{h^{\alpha }}
\left[ \sum\limits_{k=-i}^{N-i}\left( \sigma C_{i+k}^{f}+\left( 1-\sigma
\right) C_{i+k}^{f+1}\right) w_{k}^{\left( \alpha \right) }\right.   \notag
\\
&&\left. +g_{L}^{f+\frac{1}{2}}s_{L}\,_{i}^{\left( \alpha \right) }+\
g_{R}^{f+\frac{1}{2}}s_{R}\,_{i}^{\left( \alpha \right) }\right] ,
\end{eqnarray}
and for the boundary nodes $x_{0}$ and $x_{N}$:
\setlength\arraycolsep{1pt}
\begin{eqnarray}
C_{0}^{f+1} &=&{\ g_{L}^{f+\frac{1}{2}},} \\
C_{N}^{f+1} &=&{\ g_{R}^{f+\frac{1}{2}}.}
\end{eqnarray}

The method is explicit for $\sigma =1$ and partially implicit
for $0<\sigma <1$ and with $\sigma =0$ being fully implicit. In literature this
method is known as the $\sigma $-method for parabolic equations. 

Above scheme described by expressions~(37)-(39) can be written in matrix form as
\begin{equation}
\mathbf{A}\cdot \mathbf{C^{f+1}}=\mathbf{B}\text{,}
\end{equation}
where
\setlength\arraycolsep{0.5pt}
\begin{equation}
\mathbf{A}=\left[ 
\begin{array}{cccccccc}
1 & 0 & 0 & 0 & \ldots  & 0 & 0 & 0 \\ 
a_{-1} & 1+a_{0} & a_{1} & a_{2} & \ldots  & a_{N-3} & a_{N-2} & a_{N-1} \\ 
a_{-2} & a_{-1} & 1+a_{0} & a_{1} & \ldots  & a_{N-4} & a_{N-3} & a_{N-2} \\ 
a_{-3} & a_{-2} & a_{-1} & 1+a_{0} & \ldots  & a_{N-5} & a_{N-4} & a_{N-2} \\ 
a_{-4} & a_{-3} & a_{-2} & a_{-1} & \ldots  & a_{N-6} & a_{N-3} & a_{N-4} \\ 
\vdots  & \vdots  & \vdots  & \vdots  & \ddots  & \vdots  & \vdots  & \vdots \\ 
a_{-N+2} & a_{-N+3} & a_{-N+4} & a_{-N+5} & \ldots  & 1+a_{0} & a_{1} & a_{2} \\ 
a_{-N+1} & a_{-N+2} & a_{-N+3} & a_{-N+4} & \ldots  & a_{-1} & 1+a_{0} & a_{1} \\ 
0 & 0 & 0 & 0 & \ldots  & 0 & 0 & 1
\end{array}
\right] ,
\end{equation}
\begin{equation}
\mathbf{B}=\left[ 
\begin{array}{c}
{\ g_{L}^{f+\frac{1}{2}}} \\ 
b_{1} \\ 
b_{2} \\ 
b_{3} \\ 
b_{4} \\ 
\vdots  \\ 
b_{N-2} \\ 
b_{N-1} \\ 
{\ g_{R}^{f+\frac{1}{2}}}
\end{array}
\right] ,
\end{equation}
with 
\setlength\arraycolsep{1pt}
\begin{eqnarray}
a_{j} &=&\left(\sigma -1 \right) K_{\alpha }\dfrac{\Delta t}{h^{\alpha }}
w_{j}^{\left( \alpha \right) },\text{ \ \ \ for }j=-N+1,\ldots ,N-1\text{,}
\\
b_{j} &=&C_{j}^{f}+K_{\alpha }\dfrac{\Delta t}{h^{\alpha }}
\left[ g_{L}^{f+\frac{1}{2}}{s_{L}}_{j}^{\left( \alpha \right) }
+g_{R}^{f+\frac{1}{2}}{s_{R}}_{j}^{\left( \alpha \right) }\right.   \notag \\
&&\left. +\ \sigma \sum\limits_{k=-j}^{N-j}C_{i+k}^{f}w_{k}^{\left( \alpha
\right) }\right] ,\text{ \ \ \ for }j=1,\ldots ,N-1\text{.}
\end{eqnarray}
and $\mathbf{C^{f+1}}$ is the vector of unknown function's values $C$ at the
time $t^{f+1}$.

Particular case of above scheme (37) is the explicit scheme (for $\sigma =1$)
which may be simplified to 
\setlength\arraycolsep{0pt}
\begin{equation}
C_{i}^{f+1}=\left\{ {
\begin{array}{ll}
g_{L}^{f+\frac{1}{2}}, & \hspace{-1.7cm}\text{for }i=0, \\ 
 K_{\alpha }\dfrac{\Delta t}{h^{\alpha }}
\left( g_{L}^{f+\frac{1}{2}}{s_{L}}_{i}^{\left( \alpha \right) }
      +g_{R}^{f+\frac{1}{2}}{s_{R}}_{i}^{\left( \alpha \right) } \right) &  \\ 
 \hspace{0.7cm} +  \sum\limits_{k=-i}^{n-i}C_{i+k}^{f}p_{k}^{\left( \alpha \right) } \text{, } &
 \hspace{-1.7cm}\text{for }i=1,\ldots ,N-1, \\ 
g_{R}^{f+\frac{1}{2}}, & \hspace{-1.7cm}\text{for }i=N,
\end{array}
}\right. 
\end{equation}
with 
$p_{k}^{\left( \alpha \right) }$ defined by formula (26).

We can observe that boundary conditions influence to all values
of the function in every node. In opposite to the second derivative over space 
which is approximated locally, the characteristic feature of Riesz-Feller 
and other fractional derivatives is dependence on values of all domain points.
For $\alpha = 2$ and $\theta = 0$ our scheme is identically 
as wide known and used the forward difference in time and central 
difference in space scheme (FTCS)~\cite{Ames,Hoffman,Majchrzak}.

The skewness parameter $\theta$ has great significance influence on the solution.
For $\alpha \rightarrow 1^{+}$ and $\theta \rightarrow \pm 1^{+}$ one can obtain
the classical hyperbolic equation, i.e. the first order wave equation (the transport equation).
In this case our scheme tends to the known Euler's forward time and central space 
(FTCS) approximation of Eqn (1). Unfortunately this is unconditionally unstable and
therefore this is disadvantage this method.

Proposed numerical scheme makes a bridge between Gaussian and Cauchy processes.
Our scheme is also a~bridge between diffusion and transport phenomena.

\section{ Simulation results}
In this section we present results of calculation. In all presented simulations
we assumed $k_{\alpha } = 1 m^{\alpha }/s$ and the length of 1D domain $l=1 m$. 
Figure~1 shows two charts over space (one in the logarithmic scale)
with absorbing boundary $\left. C\left( x,t\right) \right\vert_{x=0}
=\left. C\left( x,t\right) \right\vert _{x=1}=0$.
On these plots solutions for different values of parameter 
$\alpha \in \left\langle \,1.01,1.5,2\right\rangle $
at time $t=0,0.01,0.3s $ for $\theta=0$ are presented.
\begin{figure}[hb] 
\begin{center}
    \includegraphics[width=0.48\textwidth]{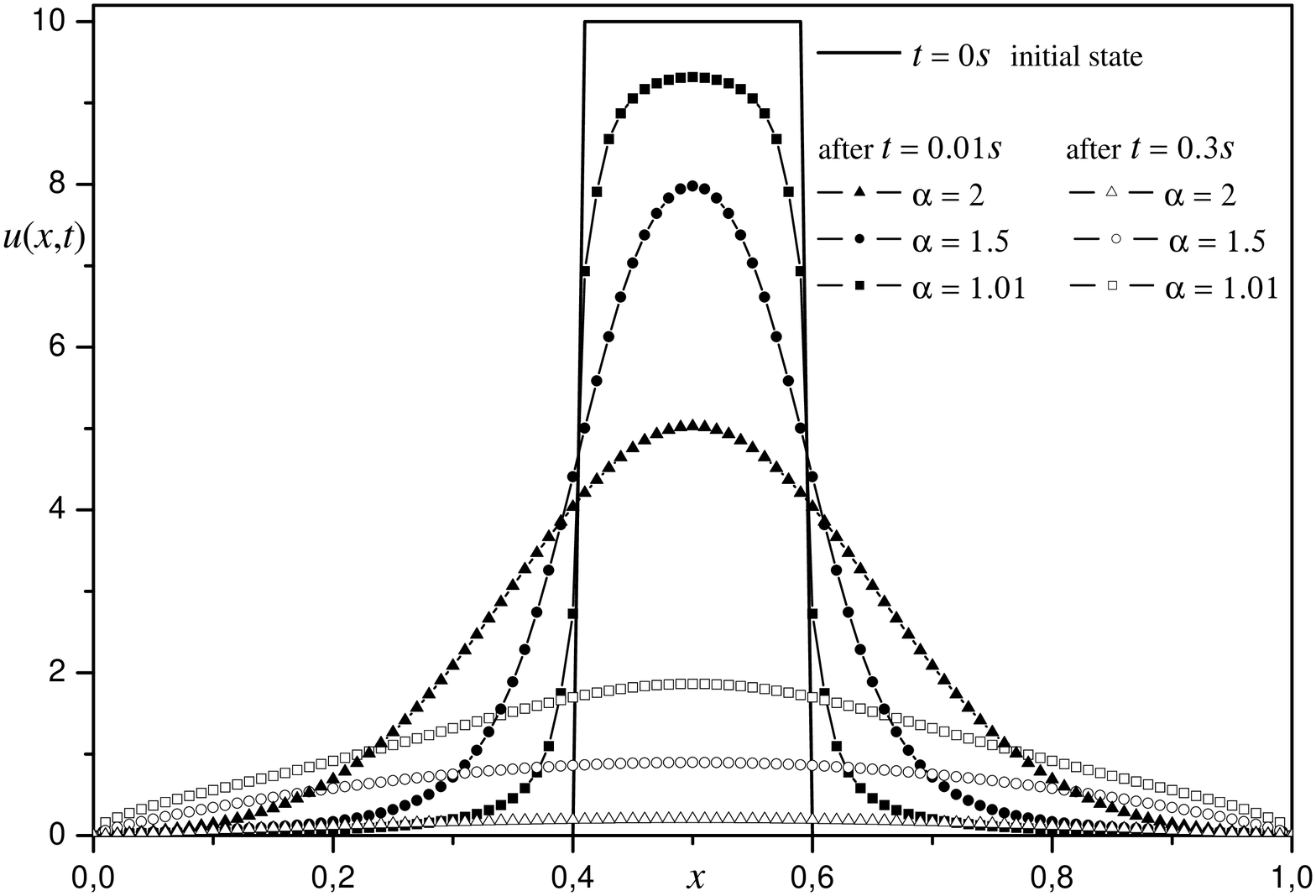}
    \includegraphics[width=0.48\textwidth]{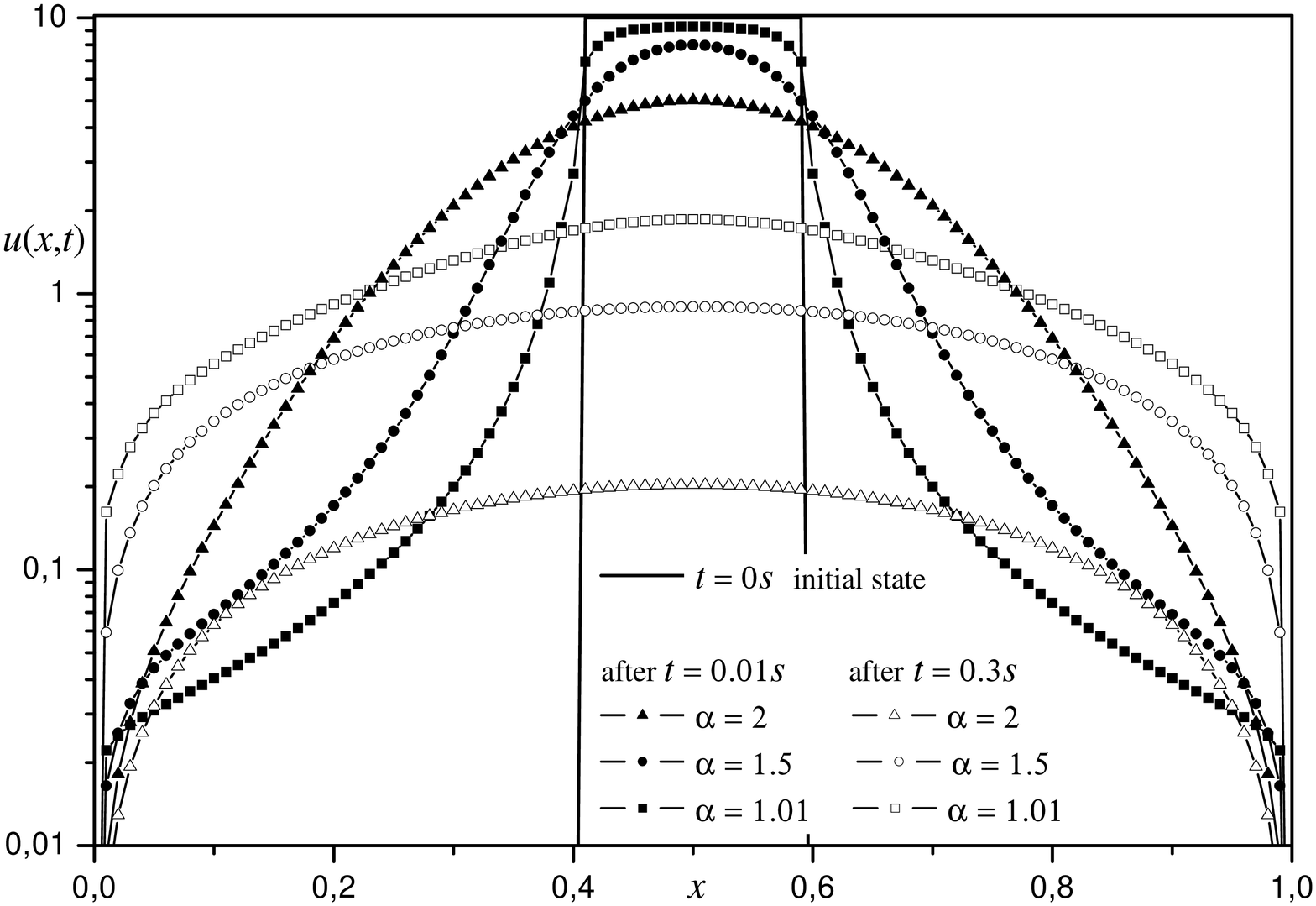}
    \caption{ Solution over space for $\alpha \in \left\langle \,1.01,1.5,2\right\rangle $ \protect\\
	a)~normal scale; \qquad b)~logarithmic scale. } 
\end{center}
\end{figure} 
Figure~2 presents another example of the solution which differs from example
presented by the Fig.~1 (boundary conditions $\left. C\left( x,t\right) \right\vert_{x=0}
=\left. C\left( x,t\right) \right\vert _{x=1}=100$ and initial condition 
$\left. C\left( x,t\right) \right\vert_{t=0} = 0$ ).
In both cases we observe diffusion process arising in different way.
The last example reflects case when the parameter of skewness is $\theta = 0.5$
and $\alpha = 1.4$. Figure~3 shows a diffusion transport process 
over space at different moments of time.
\begin{figure}[h] 
\begin{center}
\includegraphics[width=0.47\textwidth]{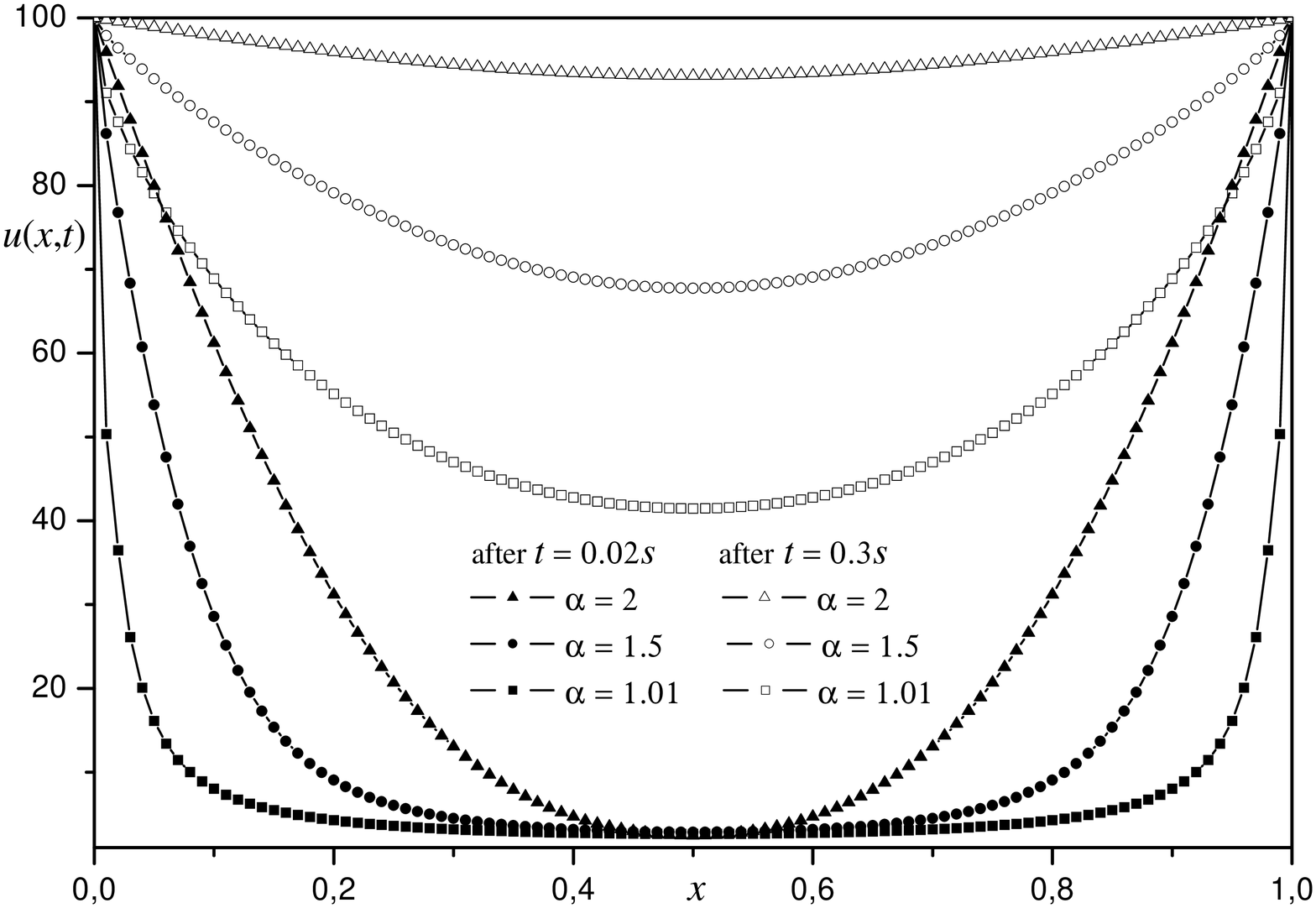}
\caption{ Solution over space for $\alpha=1.01, 1.5, 2$.  }    
\vspace{0.5cm}
 \includegraphics[width=0.47\textwidth]{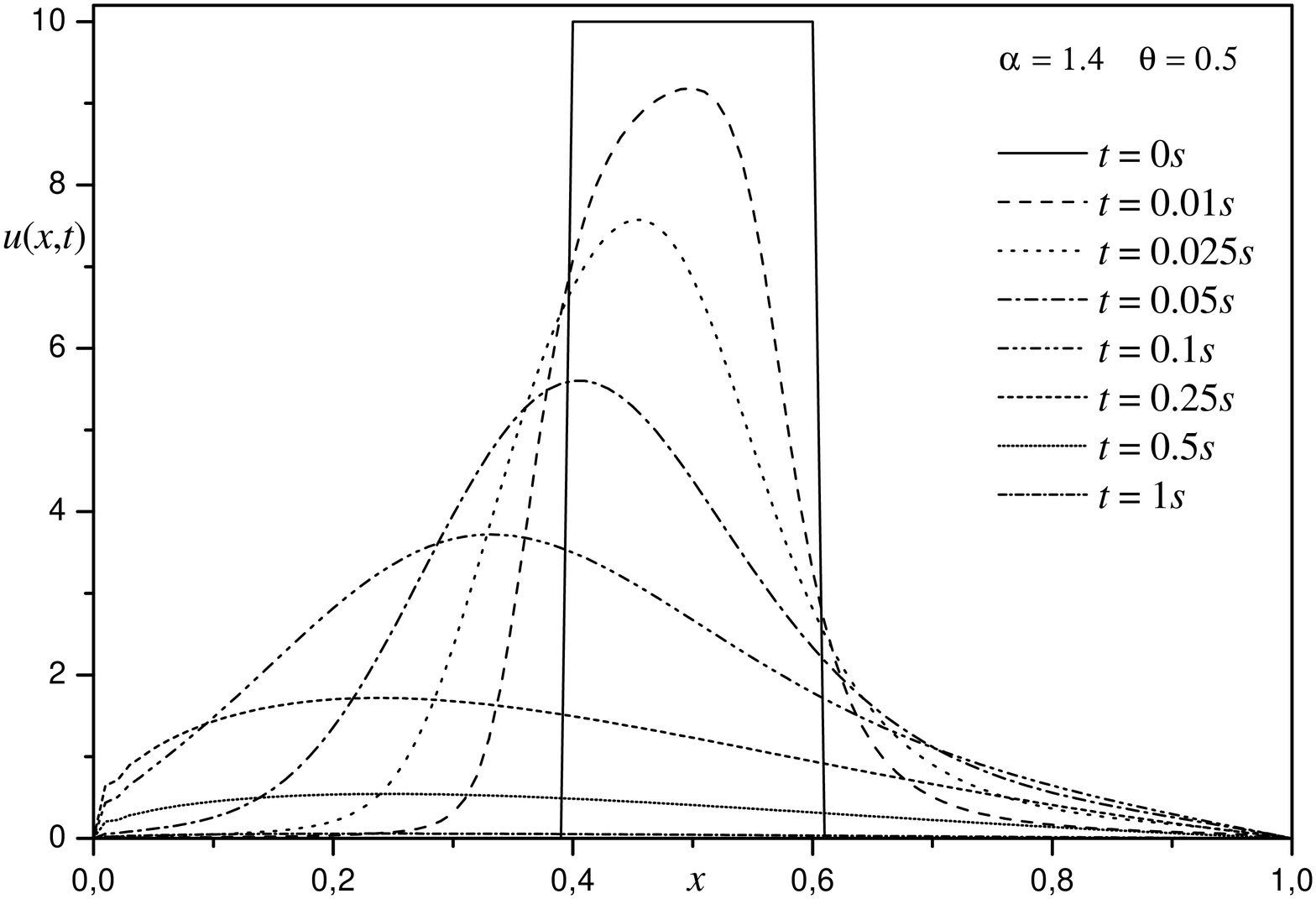}
 \caption{Solution over space for $\alpha=1.4$ and $\theta =0.5$. }
\end{center}
\end{figure}

\section{ Conclusions}
In summary, we proposed the fractional finite difference met\-hod 
for fractional diffusion equation with the Riesz-Feller fractional derivative
which is extension to the standard diffusion.
We analysed a linear case of diffusion equation and in the future 
we will work on non-linear cases. 
We obtained the implicit and explicit FDM schemes 
which generalise classical schemes of FDM for the diffusion 
equation. Moreover, for $\alpha = 2$ our solution equals to the classical
finite difference method.

Analysing plots included in this work,
we can see that in the case $\alpha < 2$ (the Levy flight)
diffusion is slower then the standard diffusion (Brownian motion) 
in the initial time. Nevertheless, when we analyse the probability density function 
we observe a long tail of distribution in the long time limit. 
In this way we can simulate same rare and extreme events which are characterised by arbitrary very large values of particle jumps.

Analysing changes in the skewness parameter $\theta$ we observed interesting
behaviour in solution.
For $\alpha \rightarrow 1^{+}$ and for $\theta \rightarrow \pm 1^{+}$
we obtained the first order wave equation. For $\theta \in (0,1)$ 
(with restrictions to order $\alpha$) we generate a class of non-symmetric 
probability density functions.

\end{document}